\documentclass[preprint,12pt]{elsarticle}
\usepackage{amsfonts,amsmath,amssymb}
\usepackage{tikz}
\usepackage[colorlinks=true, citecolor=black, linkcolor=black, urlcolor=black]{hyperref}
\usepackage{todonotes}
\newcommand{\commA}[1]{\marginpar{%
\begin{color}{red}
\vskip-\baselineskip %raise the marginpar a bit
\raggedright\footnotesize
\itshape\hrule \smallskip A: #1\par\smallskip\hrule\end{color}}}

\newtheorem{thm}{Theorem}[section]
\newtheorem{lem}[thm]{Lemma}
\newtheorem{co}[thm]{Conjecture}

\newtheorem{coex}[thm]{Counterexample}
\newdefinition{df}{Definition}[section]
\newdefinition{rem}{Remark}[section]
\newdefinition{ex}{Example}[section]
\newproof{pf}{Proof}
\newproof{pot}{Proof of Theorem}
\numberwithin{equation}{section}

\journal{
}

\begin{document}

\begin{frontmatter}

\title{The generalized characteristic polynomial, corresponding resolvent and their application
\footnote{\bf To all fearless Ukrainians defending not only their country, but the whole civilization 
against putin's 
\href{https://en.wikipedia.org/wiki/Ruscism}{rashism}
}
}
\author{A.V.~Kosyak\corref{cor1}}
%\fnref{fn1}
\ead{kosyak02@gmail.com}
\address{Institute of Mathematics, Ukrainian National Academy of Sciences,\\
3 Tereshchenkivs'ka Str., Kyiv, 01601, Ukraine}
\address{
London Institute for Mathematical Sciences,\\
21 Albemarle St, London W1S 4BS, UK
}

\cortext[cor1]{Corresponding author}

\begin{abstract}
We introduced previously 
%in \cite{Kos_B_09} 
the generalized characteristic polynomial  defined by$
P_C(\lambda)={\rm det}\,C(\lambda),$ where $C(\lambda)=C+{\rm diag}\big(\lambda_1,\dots,\lambda_n\big)$
%\sum_{k=1}^n\lambda_kE_{kk}+C.
%C(\lambda)={\rm det}\Big(C+{\rm diag}\big(\lambda_1,\dots,\lambda_n\big)\Big)
for $C\in {\rm Mat}(n,\mathbb C)$ and $\lambda=(\lambda_k)_{k=1}^n\in \mathbb C^n$ and gave the explicit formula for  $P_C(\lambda)$. 
%Lemma 3.1. 
In this article we define an analogue of the resolvent $C(\lambda)^{-1}$, calculate it and the expression $(C(\lambda)^{-1}a,a)$ for $a\in \mathbb C^n$ explicitly.
\iffalse
and  show  that 
$
1+(C(\lambda)^{-1}a,a)=
{\rm det}\big(C(\lambda)+a\otimes a\big)/{\rm det}\,C(\lambda).
%\big(C(\lambda)^{-1}a,a\big)=\frac{{\rm det}(I+\gamma(y_1,y_2,\dots,y_m))}{{\rm det}(I+\gamma(y_2,\dots,y_m))}-1,
$
%where $(y_k)_{k=1}^m$ are defined by $C,\,\lambda$ and $a\in \mathbb C^n$.
%see the proof of Theorem 5.7.
\fi
%%%
The obtained formulas and their variants were applied 
to the proof of the irreducibility of unitary representations of some infinite-dimensional groups.
%, for detail see \cite{Kos_B_09,KosJFA17,KosMor-Arx23}. 
\end{abstract}


\end{frontmatter}
%%(((14.11.12
\newpage
\tableofcontents

%%%%%%%%%%%%%%%%%%%%%%%%%%%%%%%%%%%%%%%%%%%%
%%%%%%%%%%%%

%%
\section{Summary of the key formulas}
 {\small
For the generalized characteristic polynomial 
$P_C(\lambda)\!=\!{\rm det}\,\Big(C+{\rm diag}\big(\lambda_1,\dots,\lambda_n\big)\!\Big)$  
we have (for notations see Definitions~\ref{df.min-cofact}, \ref{df.adjunct} and Remark~\ref{r.lam(alpha)}):
\begin{eqnarray*}
&&
P_C(\lambda)={\rm det}\,C(\lambda)=
\sum_{\emptyset\subseteq\alpha\subseteq\{1,2,\dots,n\}}\lambda_\alpha
A^\alpha_\alpha(C)=\Big(\prod_{k=1}^n\lambda_k\Big)\sum_{\emptyset\subseteq\alpha\subseteq\{1,2,\dots,n\}}
\frac{M^\alpha_\alpha(C)}{\lambda_\alpha},\\
&&
C(\lambda)^{-1}=\frac{1}{P_C(\lambda)}\Big(\prod_{k=1}^n\lambda_k\Big) 
\sum_{\emptyset\not=
\alpha\subseteq\{1,2,\dots,n\}}
\frac{A^T(C_\alpha)}{\lambda_\alpha},\\
&&
\big(
C(\lambda)^{-1}
a,a\big)=\frac{1}{P_C(\lambda)}\Big(\prod_{k=1}^n\lambda_k\Big) 
\sum_{\emptyset\not=
\alpha\subseteq\{1,2,\dots,n\}}\frac{\big(A^T(C_\alpha)a_\alpha,a_\alpha\big)}{\lambda_\alpha},\\
&&
1+(C(\lambda)^{-1}a,a)=\frac{{\rm det}\big(C(\lambda)+a\otimes a\big)}{{\rm det}\,C(\lambda)},\quad\text{where}\quad
a\otimes a=(a_ka_r)_{k,r=1}^n.
%\frac{{\rm det}(I+\gamma(y_1,y_2,\dots,y_m))}{{\rm det}(I+\gamma(y_2,\dots,y_m))}-1
\end{eqnarray*}
Another presentation of $1+(C(\lambda)^{-1}a,a)$, that we will use,  is given in   Theorem~\ref{t.m}.
}
\section{Characteristic polynomials}
Consider an 
$n\times n$ matrix $C$. The {\it characteristic polynomial} of 
$C$, denoted by $p_C(t)$
is the polynomial defined by 
$p_C(t)=\det(tI-C)$,
where $I$ denotes the 
$n\times n$ identity matrix.
%, see \cite{wiki-char-pol}. 
By the {\it Cayley--Hamilton theorem} we have  $p_C(C)=0$.
%, \cite{wiki-Cayl-Ham}. 
Some authors define the characteristic polynomial as
$p_C(t)=\det(C-tI)$. For a 
$2\times 2$ matrix 
$C$, the characteristic polynomial is thus given by
\begin{equation*}
p_C(t)=t^2-t\,{\rm tr}(C)+{\rm det}\,C. 
\end{equation*}
Using the language of {\it exterior algebras}, the characteristic polynomial of an 
$n\times n$ matrix  $C$ may be expressed as
\begin{equation}
\label{char-pol.1} 
 p_C(t)=\sum _{k=0}^{n}t^{n-k}(-1)^{k}{\rm tr} \left(\textstyle \bigwedge ^{k}C\right)=\sum _{k=0}^{n}t^{n-k}(-1)^{k}c_k,
\end{equation}
where 
$
c_k={\rm tr} \left(\bigwedge ^{k}C\right)$ is the trace of the 
$k^{th}$ {\it exterior power} of 
$C$, which has dimension 
$ {\binom {n}{k}}$. This trace may be computed as {\it the sum of all principal minors of $C$ of size  $k$} (see Definition~\ref{df.min-cofact} 
and Remark~\ref{r.lam(alpha)}):
\begin{equation}
c_k=\sum_{
\emptyset\subseteq
\alpha\subseteq\{1,2,\dots,n\},\,
\vert \alpha\vert=k}\lambda_\alpha
M^\alpha_\alpha(C).
\end{equation}
The recursive {\it Faddeev–LeVerrier algorithm} computes these coefficients more efficiently \cite{How98}.
%,wiki-Fad-LeVerd}.
%
When the characteristic of the field of the coefficients is 
$0$,
each such trace may alternatively be computed as a single determinant, 
that of the 
$k\times k$ matrix,
%tr
\begin{equation}
 \label{tr(ext)^k}
c_k=
{\rm tr} \left(\textstyle \bigwedge^{k}C\right)=\frac {1}{k!}
{\begin{vmatrix}
\!{\rm tr}\,C&k-1&0&\cdots &0\\ {\rm tr}\,C^{2}& {\rm tr}\, C&k-2&\cdots &0\\
\vdots &\vdots &&\ddots &\vdots \\
\,\,{\rm tr}\,C^{k-1}& {\rm tr}\,C^{k-2}&&\cdots &1\\
{\rm tr}\,C^{k}& {\rm tr}\,C^{k-1}&&\cdots &{\rm tr}\,C\end{vmatrix}}~.
\end{equation}
Theorem~\ref{t.(C^{-1}a,a)}, formula 
\eqref{C^{-1}(lam)} gives the expression for $
%C(\lambda)^{-1}=
\Big(C+{\rm diag}\big(\lambda_1,\dots,\lambda_n\big)\!\Big)^{-1}$. In particular,
for {\it resolvent} $(tI-C)^{-1}$  we have    
\begin{equation}
 \label{Res_C(t)}
(tI-C)^{-1}=\frac{1}{p_C(t)}
\Big[
\sum_{k=1}^nt^{n-k}(-1)^{k+1}
\sum_{
%\emptyset\subseteq
\alpha\subseteq\{1,2,\dots,n\},\,
\vert \alpha\vert=k}
A^T(C_\alpha)
\Big],
\end{equation}
where notation $A(C_\alpha)$ are defined in Definition~\ref{df.adjunct}. For a $3\times 3$ matrix $C$ we have for example
\begin{equation}
\label{Res_C(t).n=3}
(tI-C)^{-1}=\frac{1}{p_C(t)}
\Big[t^2\sum_{k=1}^3A^T(C_k)-t\sum_{1\leq k<r\leq 3}^3A^T(C_{kr})+A^T(C_{123})
\Big].
\end{equation}
\section{The generalized characteristic polynomial
and its properties}
\label{sec.gen.har.pol}
\index{polynomial!characteristic!generalized}
\begin{df}
\label{d.G_k(lambda)} 
For a matrix $C\in {\rm
Mat}(n,{\mathbb C})$ and $\lambda=(\lambda_k)_{k=1}^n \in {\mathbb C}^n$ define 
the {\it generalization of the characteristic polynomial} $p_C(t)={\rm det}\,(tI-C)
,\,t\in {\mathbb C}$ as follows:
\begin{equation}
\label{P_C(lambda)}
P_C(\lambda)={\rm det}\,C(\lambda),\quad\text{where}\quad
C(\lambda)=
{\rm diag}\big(\lambda_1,\dots,\lambda_n\big)+C.
\end{equation}
\end{df}
\begin{df}
 \label{df.min-cofact} 
For a matrix $C\in {\rm Mat}(n,\mathbb C),\,\,a\!\in \mathbb C^n$  and fixed $1\leq i_1<i_2<\dots i_r\leq n$ rows
and $1\leq j_1<j_2<\dots j_r\leq n$ columns $1\!\leq\! r\!\leq\! n$
denote by 
$$
M^{i_1i_2\dots i_r}_{j_1j_2\dots j_r}(C)\quad\text{ and} \quad
A^{i_1i_2\dots i_r}_{j_1j_2\dots j_r}(C)
$$
the corresponding {\it minors} and {\it cofactors} of the matrix $C$. 
\end{df}
\begin{lem}
%[\cite[Ch.1.4.3]{Kos_B_09}]
\label{l.detC-LI}
{\rm (\cite[Ch.1.4.3]{Kos_B_09})}
For the generalized characteristic polynomial $P_C(\lambda)$
of
$C\!\in\!{\rm Mat}(n,{\mathbb C})$ and
 $\lambda=(\lambda_1,\lambda_2,...,\lambda_n)\in {\mathbb
C}^n$ we have 
\begin{equation}
\label{detC-LI}
P_C(\lambda)=
{\rm det}\,C+
\sum_{r=1}^n\sum_{1\leq i_1<i_2<...<i_r\leq
n}\lambda_{i_1}\lambda_{i_2}...\lambda_{i_r}A^{i_1i_2...i_r}_{i_1i_2...i_r}(C).
\end{equation}
\end{lem}
\begin{rem}
\label{r.lam(alpha)}
If we set
$\lambda_\alpha=\lambda_{i_1}\lambda_{i_2}\cdots\lambda_{i_r}$, 
where
$\alpha=\{i_1,i_2,\dots,i_r\}$ and
$A^\alpha_\alpha(C)=A^{i_1i_2...i_r}_{i_1i_2...i_r}(C),\,\,
M^\alpha_\alpha(C)=M^{i_1i_2...i_r}_{i_1i_2...i_r}(C),\,\,
\lambda_\emptyset=1,\,\,A^\emptyset_\emptyset(C)={\rm det}\,C$ and $\vert \alpha\vert=r$,
(see Definition~\ref{df.adjunct}) we
may write (\ref{detC-LI}) as follows:
\begin{equation}
\label{A.detC-LI.2} 
%G_m(\lambda)
P_C(\lambda)={\rm det}\,C(\lambda)=
\sum_{\emptyset\subseteq\alpha\subseteq\{1,2,\dots,n\}}\lambda_\alpha
A^\alpha_\alpha(C).
\end{equation}
Writing $\widehat{\alpha}=\{1,2,\dots,n\}\setminus \alpha$, we have $A^\alpha_\alpha(C)=M^{\widehat{\alpha}}_{\widehat{\alpha}}(C)$, hence
\begin{equation}
\label{M.detC-LI.2}
%G_m(\lambda)
P_C(\lambda)={\rm det}\,C(\lambda)=
\Big(\prod_{k=1}^n\lambda_k\Big)\sum_{\emptyset\subseteq\alpha\subseteq\{1,2,\dots,n\}}
\frac{M^\alpha_\alpha(C)}{\lambda_\alpha}.
\end{equation}
\end{rem}
%
%%%%%%%%%%%%%%%%%}
\section{Gram determinants and Gram matrices}
\begin{df}
\label{d.Gram-det}
Gram determinants were introduced in 1879 by J.P. Gram \cite{Gram1879}.   
For vectors $x_1,x_2,\dots,x_n$ in some
Hilbert space $H$ the {\it Gram
matrix} $\gamma(x_1,x_2,\dots, x_n)$ is defined  by the formula (see also \cite{Gan58}, Chap IX, \S 5)
$$
\gamma(x_1,x_2,\dots, x_n)=(x_k,x_m)_{k,m=1}^n.
$$
The determinant of this matrix is called the {\it Gram  determinant} for
the vectors $x_1,x_2,..., x_n$ and is denoted by
	%$G(x_1,x_2,..., x_m)$.
$\Gamma(x_1,x_2,\dots, x_n)$ 	
\begin{equation}
\label{Gram-det}
\Gamma(x_1,x_2,\dots, x_n):={\rm det}\,\gamma(x_1,x_2,\dots, x_n).
\end{equation}
Some authors  use the notation $G(x_1,x_2,\dots, x_n)$.
\end{df}
\index{determinant!Gram}
\index{matrix!Gram}
\index{parallelotope}
\begin{rem}
\label{Vol.par-pe} 
A Gram determinant is equal to the square of the $n-$dimensional volume of the {\it parallelotope} constructed on $x_1,x_2,...,x_n$.
\end{rem}
\index{Gram determinant}
\index{Gram matrix}
\index{parallelotope}

Fix some notations
\begin{equation}
\label{X(mn)}
X=X_{mn} =\left(
\begin{array}{cccc}
x_{11}&x_{12}&...&x_{1n}\\
x_{21}&x_{22}&...&x_{2n}\\
...   &...   &...&...\\
x_{m1}&x_{m2}&...&x_{mn}
\end{array}
\right),
\end{equation}
%Set 
%\commA{corrected $y_r$, 10.04.20}
\begin{equation}
\label{x_k,y_r=}
x_k=(x_{1k},x_{2k},\dots,x_{mk})\in{\mathbb R}^m,
%\,1\leq k\leq m,
\quad
y_r=(x_{r1},x_{r2},\dots,x_{rn})\in{\mathbb R}^n.
%\,1\leq r\leq m,
\end{equation}
Then, obviously, we get
\begin{equation}
\label{X^*X}
X^*X=
 \left(
\begin{array}{cccc}
(x_1,x_1)&(x_1,x_2)&\dots&(x_1,x_n)\\
(x_2,x_1)&(x_2,x_2)&\dots&(x_2,x_n)\\
\dots   &\dots   &\dots&\dots\\
(x_n,x_1)&(x_n,x_2)&\dots&(x_n,x_n)
\end{array}
\right)=\gamma(x_1,x_2,\dots,x_n),
\end{equation}
\begin{equation}
\label{XX^*}
XX^*= \left(
\begin{array}{cccc}
(y_1,y_1)&(y_1,y_2)&\dots&(y_1,y_m)\\
(y_2,y_1)&(y_2,y_2)&\dots&(y_2,y_m)\\
\dots   &\dots   &\dots&\dots\\
(y_m,y_1)&(y_m,y_2)&\dots&(y_m,y_m)
\end{array}
\right)=\gamma(y_1,y_2,..., y_m).
\end{equation}
Therefore, we have
\begin{equation}
 \label{d(X^*X)=d(XX^*)}
\Gamma(x_1,x_2,..., x_n)={\rm det}(X^*X)={\rm det}(XX^*)=\Gamma(y_1,y_2,..., y_m).
\end{equation}
%%%%%
%
%https://www.encyclopediaofmath.org/index.php/Gram_determinant}
%
\subsection{How far is  a vector from a hyperplane}
\label{s.dist}
We start with a classical result, see, e.g.,
% Gantmacher 
\cite{Gan58}.
Consider the hyperplane $V_n$ generated by $n$ arbitrary vectors $f_1,\dots,f_{n}$ in some Hilbert space  $H$. 
\begin{lem} [\cite{Gan58,AhiGlaz93}]
\label{l.d(f,v_n)}
The square of the distance $d(f_0,V_n)$ of a vector $f_0$ from the hyperplane $V_n$ is given by the ratio of two Gram determinants (see Definition \ref{d.Gram-det})
\begin{equation}
\label{d(f,v_n)}
d^2(f_0,V_n)=
%\frac{G(f_0,f_1,f_2,\dots, f_n,)}{G(f_1,f_2,\dots, f_n)}.
\frac{\Gamma(f_0,f_1,f_2,\dots,f_n)}{\Gamma(f_1,f_2,\dots, f_n)}.
\end{equation}
\end{lem}
%
%\commA{11.10.23, to give the proof?}
%
\begin{pf} 
%For the convenience of reader we give here the proof,
We 
follow closely the book by Axiezer and Glazman \cite{AhiGlaz93}.
Set $f=\sum_{k=1}^nt_kf_k\in V_n$ and $h=f-f_0$. Since $h$ should be orthogonal to $V_n$ we conclude that
$f_r\perp h$,  i.e., $(f_r,h)=0$ for all $r$, or
\begin{equation}
 \label{f_r-perp-h}
\sum_{k=1}^nt_k(f_r,f_k)=(f_r,f_0),\quad 1\leq r\leq n. 
\end{equation}
Set $A=\gamma(f_1,f_2,\dots,f_n)$ and $b=(f_k,f_0)_{k=1}^n\in \mathbb R^{n} $. By definition we have 
\begin{equation}
\label{d^2} 
d^2=\min_{f\in V_n}\Vert f-f_0\Vert^2=(At,t)-2(t,b)+(f_0,f_0).
\end{equation}
Since $d^2=(h,h)=(f_0,h)$ we conclude that $d^2=\sum_{k=1}^nt_k(f_0,f_k)-(f_0,f_0)$ or
\begin{equation}
 \label{d^2=(f_0,h)}
\sum_{k=1}^nt_k(f_0,f_k)=(f_0,f_0)-d^2. 
\end{equation}
So we have the system of equations:
\begin{equation}
\label{syst}
\left\{
\begin{array}{ccc}
t_1(f_1,f_1)+t_2(f_1,f_2)+\dots +t_n(f_1,f_n)&=&(f_1,f_0)\\
t_1(f_2,f_1)+t_2(f_2,f_2)+\dots +t_n(f_2,f_n)&=&(f_2,f_0)\\
\dots&&\\
t_1(f_n,f_1)+t_2(f_n,f_2)+\dots +t_n(f_n,f_n)&=&(f_n,f_0)\\
t_1(f_0,f_1)+t_2(f_0,f_2)+\dots +t_n(f_0,f_n)&=&(f_0,f_0)-d^2
\end{array}
\right..
\end{equation}
Excluding $t_k$ from the system we get $d^2=\frac{\Gamma(f_0,f_1,f_2,\dots, f_n,)}{\Gamma(f_1,f_2,\dots, f_n)}$. 
Formula \eqref{d(f,v_n)} also follows from Remark~\ref{Vol.par-pe}.
\qed\end{pf}
\begin{rem}
\label{r.At=b}
From the system \eqref{syst} we conclude that $At=b$, where $b=(f_k,f_0)_{k=1}^n\in \mathbb R^{n} $, hence $t=A^{-1}b$. 
By \eqref{d^2} we get 
\begin{equation}
\label{At=b} 
d^2=(f_0,f_0)-(A^{-1}b,b)=\frac{\Gamma(f_0,f_1,f_2,\dots, f_n)}{\Gamma(f_1,f_2,\dots, f_n)}.
\end{equation}
See also \cite[Chap. 4.3, Lemma 4.3.2]{Kos_B_09}.
\end{rem}

\section{The explicit expression for $C^{-1}(\lambda)$ and  $(C^{-1}(\lambda)a,a)$}
\label{sec.C^{-1}(lam)}
Fix $C\in {\rm Mat}(n,\mathbb C),$ $a\in \mathbb C^n$ and $\lambda\in \mathbb C^n$. 
Our aim is to find the explicit formulas for  
%${\rm det}\,C(\lambda)$,
$C(\lambda)^{-1}$ and  $(C(\lambda)^{-1}a,a)$, where $C(\lambda)$ is defined by \eqref{P_C(lambda)}.
\begin{df}
\label{df.adjunct} 
For $\alpha=\{i_1,i_2,\dots, i_r\}\subset\{1,2,\dots,n\}$ set $M(\alpha)(C)=
M^\alpha_\alpha(C)$.
%M(i_1i_2\dots i_r)(C)=M^{i_1i_2\dots i_r}_{i_1i_2\dots i_r}(C)
Let also $C_\alpha=C_{i_1i_2\dots i_r}$ be the corresponding {\it submatrix} of the matrix $C$ and 
$a_\alpha
%a_{i_1i_2\dots i_r}
=(a_{i_1},a_{i_2},\dots,a_{i_r})$.
\index{$C_{i_1i_2\dots i_r}$, submatrix of the matrix $C$}
The elements of the matrix $C_\alpha$ are
on the intersection of $i_1,i_2,\dots, i_r$ rows and column of the matrix $C$.
Denote by $A(C_{i_1i_2\dots i_r})$ the matrix of the cofactors of the first order of the matrix $C_{i_1i_2\dots i_r}$, 
another name is {\it adjugate matrix}, occasionally known as {\it adjunct matrix}:
\begin{equation}
\label{A(C)}
A(C_{i_1i_2\dots i_r})=(A^{i}_j(C_{i_1i_2\dots i_r}))_{1\leq i,j\leq r}.
%\quad{set}\quad A(C_k)=1\,\,\text{for}\,\,1\leq k\leq n.
\end{equation}
The minor of order zero
is often defined to be 1, and therefore, we set  $A(C_k)=1$ for $1\leq k\leq n$.
As usual, denote by $B^T$ the {\it matrix transposed} to $B$.
\end{df}
\index{matrix!transposed}
\index{matrix!adjugate}
\index{matrix!adjanct}

Let $n=3$, then $A(C_{123})=A(C)$ is the following matrix: 
\begin{equation}
 \label{A(C)3}
A(C)=A(C_{123})=
\left(
\begin{array}{ccc}
A^1_1&A^1_2&A^1_3\\
A^2_1&A^2_2&A^2_3\\
A^3_1&A^3_2&A^3_3
\end{array}
\right)=
\left(
\begin{array}{ccc}
M^{23}_{23}&-M^{23}_{13}&M^{23}_{12}\\
-M^{13}_{23}&M^{13}_{13}&-M^{13}_{12}\\
M^{12}_{23}&-M^{12}_{13}&M^{12}_{12}
\end{array}
\right),
\end{equation}
we write $M^{ij}_{rs}$ instead of $M^{ij}_{rs}(C)$ and $A^i_j$ instead of $A^i_j(C)$.
\begin{rem}
\label{r.C_{12}=}
If ${\rm det}\,C_{i_1i_2\dots i_r}\not=0$ we have
\begin{equation}
 \label{A(C_{12})=}
A^T(C_{i_1i_2\dots i_r})={\rm det}\,C_{i_1i_2\dots i_r}\Big(C_{i_1i_2\dots i_r}\Big)^{-1}.
\end{equation}
\end{rem}
In what follows we need to consider the submatrix $A^T
\big(C_{i_1i_2\dots i_r}\big),\,\,1\leq r\leq n,$ of the matrix $C\in {\rm Mat}(n,\mathbb C)$ as an 
{\it appropriate element of} ${\rm Mat}(n,\mathbb C)$. 
\begin{thm}
%[\cite{Kos-hpl-arx23}]
\label{t.(C^{-1}a,a)}
For the matrix $C(\lambda)$ defined by 
 %\eqref{C(lam)}
\eqref{P_C(lambda)},  $a\in \mathbb C^n$  and $\lambda\in \mathbb C^n$ we have 
 {\small
 %
%   \begin{equation}
% \label{Del_n(lam,C)}
% P_C(\lambda)=\Big(\prod_{k=1}^n\lambda_k\Big)\sum_{r=1}^n\sum_{1\leq i_1<i_2<\dots <i_r\leq n}
%\frac{M(i_1i_2\dots i_r)}{\lambda_{i_1}\lambda_{i_2}\dots \lambda_{ i_r}},
%\end{equation}
 \begin{equation}
 \label{C^{-1}(lam)}
C(\lambda)^{-1}
=\frac{1}{P_C(\lambda)}\Big(\prod_{k=1}^n\lambda_k\Big)
\sum_{r=1}^n
\sum_{1\leq i_1<i_2<\dots i_r\leq n}
\frac{A^T(C_{i_1i_2\dots i_r})}{\lambda_{i_1}\lambda_{i_2}\dots \lambda_{ i_r}},
 \end{equation} 
 \begin{equation}
 \label{(C^{-1}(lam)a,a)}
\big(
C(\lambda)^{-1}
a,a\big)=\frac{1}{P_C(\lambda)}\Big(\prod_{k=1}^n\lambda_k\Big) 
\sum_{r=1}^n
\sum_{1\leq i_1<i_2<\dots <i_r\leq n}\!\!\!\!\!\!
\frac{(A^T(C_{i_1i_2\dots i_r})a_{i_1i_2\dots i_r},a_{i_1i_2\dots i_r})}{\lambda_{i_1}\lambda_{i_2}\dots \lambda_{ i_r}},
 \end{equation}
 \begin{equation}
 \label{(C^{-1}(lam)a,a).1}
\big(
C(\lambda)^{-1}
a,a\big)=\frac{1}{P_C(\lambda)}\Big(\prod_{k=1}^n\lambda_k\Big) 
\sum_{
%\emptyset\subseteq
\alpha\subseteq\{1,2,\dots,n\},\,\vert \alpha\vert\geq 1}\frac{\big(A^T(C_\alpha)a_\alpha,a_\alpha\big)}{\lambda_\alpha}.
 \end{equation} 
} 
\end{thm}
\begin{pf} 
%The first statement follows from \eqref{M.detC-LI.2}.
For $n=2$ we have
{\small
\begin{equation}
 \label{A(C)2}
C(\lambda)\!=\!
 \left(\!\!
\begin{array}{cc}
c_{11}+\lambda_1&c_{12}\\
c_{21}&c_{22}+\lambda_2
\end{array} 
\!\!\right),\,\,\,A^T(C_{12})\!=\!
 \left(\!\!
\begin{array}{cc}
A^1_1&A^2_1\\
A^1_2&A^2_2
\end{array}
\!\!\right)\!=\!
\left(\!\!
\begin{array}{cc}
c_{22}&-c_{12}\\
-c_{21}&c_{11}
\end{array}
\!\!\right),
\end{equation}
\begin{eqnarray}
\nonumber
&&
C(\lambda)^{-1}
\!=\!\frac{1}{P_C(\lambda)}
 \left(\!
\begin{array}{cc}
c_{22}+\lambda_2&-c_{12}\\
-c_{21}&c_{11}+\lambda_1
\end{array} 
\!\right)
\!=\!\frac{1}{P_C(\lambda)}\left[
\left(\!
\begin{array}{cc}
\lambda_2&0\\
0&\lambda_1
\end{array}
\!\right)
\!+\!A^T(C_{12})\right]=\\
\label{C^{-1}(la).2}
&&
\frac{\lambda_1\lambda_2}{P_C(\lambda)}
\left[\left(\!\begin{array}{cc}\lambda_1^{-1}&0\\
0&\lambda_2^{-1}\end{array}\!\right)\!+\!\frac{A^T(C_{12})}{\lambda_1\lambda_2}\right]\!
=\!\frac{\lambda_1\lambda_2}{P_C(\lambda)}
\left[
\sum_{k=1}^2
\frac{A^T(C_k)}{\lambda_k}
\!+\!\frac{A^T(C_{12})}{\lambda_1\lambda_2}\right],\,\,
\end{eqnarray}
}
%
%\commA{\eqref{C^{-1}(la).2}}

recall that $A^T(C_1)=
\left(\begin{smallmatrix}
1&0\\
0&0
\end{smallmatrix}\right),\,\,
A^T(C_2)=
\left(\begin{smallmatrix}
0&0\\
0&1
\end{smallmatrix}\right)
$.
Therefore,
%
%\commA{\eqref{(C^{-1}(la)a,a).2}}
%
\begin{eqnarray}
\nonumber
 &&\big(
C(\lambda)^{-1}
a,a\big)=\frac{1}{P_C(\lambda)}\left[(c_{22}+\lambda_2)a_1^2-(c_{12}+c_{21})a_1a_2+(c_{11}+\lambda_1)a_2^2
 \right]\\
 \nonumber
 &&=\frac{1}{P_C(\lambda)}\left[\lambda_2a_1^2+\lambda_1a_2^2
+ c_{22}a_1^2+c_{11}a_2^2-(c_{12}+c_{21})a_1a_2\right]=\\
 \label{(C^{-1}(la)a,a).2}
&&\left(1\!+\!\frac{M(1)}{\lambda_1}+\frac{M(2)}{\lambda_2}+\frac{M(12)}{\lambda_1\lambda_2}\right)^{-1}\!\!\left[\!
\frac{a_1^2}{\lambda_1}\!+\!\frac{a_2^2}{\lambda_2}\!+\frac{(A^T(C_{12})a_{12},a_{12})}{\lambda_1\lambda_2}
 \!\right].
\end{eqnarray}
For $n=3$ we have by \eqref{M.detC-LI.2}
%
%\commA{\ref{D_3(la,C)}}
%
\begin{eqnarray}
\label{D_3(la,C)}
 &&P_C(\lambda)\!=\!\lambda_1\lambda_2\lambda_3
 \Big(1\!+\!\sum_{k=1}^3\frac{M(k)}{\lambda_k}+\sum_{1\leq k<r\leq 3}\frac{M(kr)}{\lambda_k\lambda_r} +
 \frac{M(123)}{\lambda_1\lambda_2\lambda_3}\Big),\\
\nonumber
&&C(\lambda)=
 \left(
\begin{array}{ccc}
c_{11}+\lambda_1&c_{12}&c_{13}\\
c_{21}&c_{22}+\lambda_2&c_{23}\\
c_{31}&c_{32}&c_{33}+\lambda_3
\end{array}
\right),\quad 
C(\lambda)^{-1}
=
\frac{1}{P_C(\lambda)}\times\\
%\end{eqnarray}\begin{eqnarray}
%
\nonumber
&& 
\left(\begin{smallmatrix}
\lambda_2\lambda_3\big(1+\frac{M^2_2}{\lambda_2}+\frac{M^3_3}{\lambda_3}+\frac{M^{23}_{23}}{\lambda_2\lambda_3}\big)&
-M^{13}_{23}-\lambda_3M^1_2&
M^{12}_{23}-\lambda_2M^1_3\\
-M^{23}_{13}-\lambda_3M^2_1&
\lambda_1\lambda_3\big(1+\frac{M^1_1}{\lambda_1}+\frac{M^3_3}{\lambda_3}+\frac{M^{13}_{13}}{\lambda_1\lambda_3}\big)&
-M^{12}_{13}-\lambda_1M^2_3\\
-M^{23}_{12}-\lambda_2M^3_1&
-M^{13}_{12}-\lambda_1M^3_2&
\lambda_1\lambda_2\big(1+\frac{M^1_1}{\lambda_1}+\frac{M^2_2}{\lambda_2}+\frac{M^{12}_{12}}{\lambda_1\lambda_2}
\big)
%\end{array}\right)
\end{smallmatrix}\right)=
\\
\nonumber
&&
\frac{\lambda_1\lambda_2\lambda_3}{P_C(\lambda)}
%\\
%
%\nonumber
%&&\times
\left(\begin{smallmatrix}
\frac{1}{\lambda_1}+\frac{M^2_2}{\lambda_1\lambda_2}+\frac{M^3_3}{\lambda_1\lambda_3}+\frac{M^{23}_{23}}{\lambda_1\lambda_2\lambda_3}&
-\frac{M^{13}_{23}}{\lambda_1\lambda_2\lambda_3}-\frac{M^1_2}{\lambda_1\lambda_2}&
\frac{M^{12}_{23}}{\lambda_1\lambda_2\lambda_3}-\frac{M^1_3}{\lambda_1\lambda_3}\\
-\frac{M^{23}_{13}}{\lambda_1\lambda_2\lambda_3}-\frac{M^2_1}{\lambda_1\lambda_2}&
\frac{1}{\lambda_2}+\frac{M^1_1}{\lambda_1\lambda_2}+\frac{M^3_3}{\lambda_2\lambda_3}+\frac{M^{13}_{13}}{\lambda_1\lambda_2\lambda_3}&
-\frac{M^{12}_{13}}{\lambda_1\lambda_2\lambda_3}-\frac{M^2_3}{\lambda_2\lambda_3}\\
-\frac{M^{23}_{12}}{\lambda_1\lambda_2\lambda_3}-\frac{M^3_1}{\lambda_1\lambda_3}&
-\frac{M^{13}_{12}}{\lambda_1\lambda_2\lambda_3}-\frac{M^3_2}{\lambda_2\lambda_3}&
\frac{1}{\lambda_3}+\frac{M^1_1}{\lambda_1\lambda_3}+\frac{M^2_2}{\lambda_2\lambda_3}+\frac{M^{12}_{12}}{\lambda_1\lambda_2\lambda_3}
\end{smallmatrix}\right).
\end{eqnarray}
Finally, we get
\begin{equation}
C(\lambda)^{-1}
\!=\!\frac{\lambda_1\lambda_2\lambda_3}{P_C(\lambda)}\left[
\sum_{k=1}^3
\frac{A^T(C_k)}{\lambda_k}
%\left(\begin{smallmatrix}
%       \lambda_1^{-1}&0&0\\
%       0&\lambda_2^{-1}&0\\
%       0&0&\lambda_3^{-1}
%      \end{smallmatrix}
%      \right)
+\!\!\!\!\sum_{1\leq r<s\leq 3}\frac{A^T(C_{rs})}{\lambda_r\lambda_s}
+\frac{A^T(C_{123})}{\lambda_1\lambda_2\lambda_3}
\right],
\end{equation}
we use \eqref{A(C)2} and \eqref{A(C)3}. Therefore,
\begin{eqnarray}
\nonumber
&& \big(
C(\lambda)^{-1}
a,a\big)=\frac{\lambda_1\lambda_2\lambda_3}{P_C(\lambda)}
\Big[ \frac{a_1^2}{\lambda_1}+\frac{a_2^2}{\lambda_2}+\frac{a_3^2}{\lambda_3}+\frac{(A^T(C_{12})a_{12},a_{12})}{\lambda_1\lambda_2}+
\\
\label{(C_3^{-1}a,a).n=3}
&&
%\!\!+
\frac{(A^T(C_{13})a_{13},a_{13})}{\lambda_1\lambda_3}\!+\!
\frac{(A^T(C_{23})a_{23},a_{23})}{\lambda_2\lambda_3}\!+\!\frac{(A^T(C_{123})a_{123},a_{123})}{\lambda_1\lambda_2\lambda_3}
\Big].\,\,
\end{eqnarray}
%
%In the general case we have by \eqref{A.detC-LI.2} and \eqref{M.detC-LI.2}
%\begin{equation*}
%\label{A.detC-LI.2} 
%\end{equation*}
%
%\commA{to prove in the general case \eqref{C^{-1}(lam)}}
%
%
%{\color{red} The general case}.
For $n=4$ we have
$$
C(\lambda)=
\left(\begin{array}{cccc}
c_{11}+\lambda_1&c_{12}&c_{13}&c_{14}\\
c_{21}&c_{22}+\lambda_2&c_{23}&c_{24}\\
c_{31}&c_{32}&c_{33}+\lambda_3&c_{34}\\
c_{41}&c_{42}&c_{43}&c_{44}+\lambda_4
\end{array}
\right).
$$
The general formulas are as follows 
%(the first formula follows from Lemma~\ref{l.detC-LI})
%
%\commA{small(}
%
 {\small
\begin{eqnarray*}
%
%P_C(\lambda)=\Big(\prod_{k=1}^n\lambda_k\Big)\sum_{r=1}^n\sum_{1\leq i_1<i_2<\dots <i_r\leq n}\frac{M(i_1i_2\dots i_r)}{\lambda_{i_1}\lambda_{i_2}\dots \lambda_{ i_r}},\\
%
&&
C(\lambda)^{-1}=\frac{1}{P_C(\lambda)}\Big(\prod_{k=1}^n\lambda_k\Big) 
\sum_{
%\emptyset\subseteq
\alpha\subseteq\{1,2,\dots,n\},\,\vert \alpha\vert\geq 1}\frac{A^T(C_\alpha)}{\lambda_\alpha},\\
%%%
%\frac{1}{P_C(\lambda)}\Big(\prod_{k=1}^n\lambda_k\Big)%\Delta_n^{-1}(\lambda,C) 
%\sum_{r=1}^n\sum_{1\leq i_1<i_2<\dots i_r\leq n}
%\frac{A^T(C_{i_1i_2\dots i_r})}{\lambda_{i_1}\lambda_{i_2}\dots \lambda_{ i_r}},\\
%%%%
&&
\big(
C(\lambda)^{-1}
a,a\big)=\frac{1}{P_C(\lambda)}\Big(\prod_{k=1}^n\lambda_k\Big) 
\sum_{
%\emptyset\subseteq
\alpha\subseteq\{1,2,\dots,n\},\,\vert \alpha\vert\geq 1}\frac{\big(A^T(C_\alpha)a_\alpha,a_\alpha\big)}{\lambda_\alpha}.
%%%
%\big(C(\lambda)^{-1}a,a\big)=\frac{1}{P_C(\lambda)}\Big(\prod_{k=1}^n\lambda_k\Big)
%\sum_{r=1}^n\sum_{1\leq i_1<i_2<\dots i_r\leq n}\!\!\!\!\frac{(A^T(C_{i_1i_2\dots i_r})a_{i_1i_2\dots i_r},a_{i_1i_2\dots i_r})}{\lambda_{i_1}\lambda_{i_2}\dots \lambda_{ i_r}},
%
\end{eqnarray*}
}
that proves 
\eqref{C^{-1}(lam)}--\eqref{(C^{-1}(lam)a,a).1}.
%{Del_n(lam,C)}
We make convention in \eqref{(C^{-1}(lam)a,a)}, that $A(C_k)\!=\!1$.
\qed\end{pf}
\begin{ex}
For the matrix  $C(\lambda)$ we have by \eqref{M.detC-LI.2}
\begin{eqnarray}
%\begin{equation}
&&
\label{D(lambda)}
C(\lambda)
= \left(
\begin{array}{cccc}
	1+\lambda_1&1          &...&1\\
	1          &1+\lambda_2&...&1\\
	&&...&\\
	1          &          1&...&1+\lambda_n
	\end{array}
	\right),\\
%	\end{equation}
\label{det(C).1}
&& 
{\rm det}\,C(\lambda)=\Big(\prod_{k=1}^n\lambda_k\Big)
\Big(1+\sum_{k=1}^n\frac{1}{\lambda_k}\Big),\\
\nonumber
&&
C(\lambda)^{-1}
=\left(1+\sum_{k=1}^n\frac{1}{\lambda_k}\right)^{-1}
\Big[
\sum_{k=1}^n\frac{A^T(C_k)}{\lambda_k}
+\sum_{1\leq k<r\leq n}\frac{A^T(C_{kr})}{\lambda_k\lambda_r}\Big],
\end{eqnarray}
where $A^T(C_{kr})=\left(\begin{smallmatrix}	                         1&-1\\
-1&1	                         \end{smallmatrix}\right)$
and $A^T(C_{krs})=0$ for $1\leq k<r<s\leq n$.
%
%\commA{10.10.23, if-fi 1145-70}
%
\end{ex}

\subsection{
The case where
%Particular case of 
$C$ is the Gram matrix}

Fix two natural numbers $n,m\in \mathbb N$ with $m\leq n$, two matrices $A_{mn}$ and $X_{mn}$, vectors $g_k\in\mathbb C^{m-1},\,\,1\leq k\leq n$
and $a\in \mathbb C^n$ as follows
\begin{equation} 
\label{A(mn)}
A_{mn}\!=\!\left(
	\begin{array}{cccc}
	a_{11}&a_{12}          &...&a_{1n}\\
	a_{21}&a_{22}          &...&a_{2n}\\
	&&...&\\
	a_{m1}&a_{m2}          &...&a_{mn}
	\end{array}
	\right),\,\, 
g_k=\left(
	\begin{array}{c}
	a_{2k}\\
	a_{3k}\\
	... \\
	a_{mk}
	\end{array}
	\right)\in \mathbb C^{m-1},\,\,\, a=(a_{1k})_{k=1}^n\in \mathbb C^n.	
\end{equation}
%
%\commA{ ref $A(3n)$ to  $A(mn)$. 15.04.21}
%
Set 
\begin{equation}
\label{C=Gram}
C=\gamma(g_1,g_2,\dots,g_n)
=\left(
 %\begin{smallmatrix}
\begin{array}{cccc}
 (g_1,g_1)& (g_1,g_2)&\dots& (g_1,g_n)\\
(g_2,g_1)& (g_2,g_2)&\dots& (g_2,g_n)\\
&&\dots&\\
(g_n,g_1)& (g_n,g_2)&\dots& (g_n,g_n)
%\end{smallmatrix}
 \end{array}
 \right).     
    \end{equation}
 We calculate $P_C(\lambda),\,\,C^{-1}(\lambda)$ and $(C^{-1}(\lambda)a,a)$ for an arbitrary $n$. 
Consider the matrix
%, its columns we denote by $x_k$ and its rows by $y_r$:
%
%\commA{\eqref{X(3n)}, \eqref{x_k,y_r}}
%
\begin{eqnarray}
 \label{X(3n)}
&&
 X_{mn}\!=\!\left(\!\!
\begin{array}{cccc}
x_{11}&x_{12}&...&x_{1n}\\
x_{21}&x_{22}&...&x_{2n}\\
&&...&\\
x_{m1}&x_{m2}&...&x_{mn}
\end{array}\!\!
\right),\quad 
\text{where}\quad
x_{rk}\!
=\!\frac{a_{rk}}{\sqrt \lambda_k},
%\,\, x_k\!=\!(x_{rk})_{r=1}^m\in \mathbb R^m,\,\,
\\
\label{x_k,y_r}
&&
\bar{x}_k\!=\!(x_{rk})_{r=2}^m=\frac{g_k}{\sqrt{\lambda_k}}\!\in\! \mathbb C^{m-1}.
%\quad y_r\!=\!(x_{rk})_{k=1}^n\!\in\! \mathbb R^n.
\end{eqnarray}
For $k\in \mathbb N$ define $\Delta(y_1,y_2,\dots,y_k)$ as follows:
\begin{equation}
\label{Delta(f_k)k}
\Delta(y_1,y_2,\dots,y_k)=
\frac{{\rm det}(I+\gamma(y_1,y_2,\dots,y_k))}{{\rm det}(I+\gamma(y_2,\dots,y_k))}-1.
\end{equation}

%\subsection{Arbitrary  $m\in\mathbb N$}
\begin{lem}
%[\cite{KorBogIze93}, Ch. XV]
%[O.~Gamayun]
\label{l.(A^{-1}a,a)}
For  $A\in {\rm GL}(n,{\mathbb C})$ and $a\in \mathbb C^n$ we have
\begin{equation}
\label{(A^{-1}a,a)} 
1+(A^{-1}a,a)=\frac{{\rm det}\big(A+a\otimes a\big)}{{\rm det}\,(A)}.
\end{equation}
\end{lem}
\begin{pf}
Define $a\otimes a$ as %$a\otimes a=$
$(a_ka_r)_{k,r=1}^n\in {\rm Mat}(n,{\mathbb C})$. Then by 
 \eqref{char-pol.1} we have 
 %
%\commA{13.10.23, to understnd the formula 
%\eqref{(C^{-1}(lam)a,a)}}
%
\begin{eqnarray*}
&&
{\rm det}\big(A+a\otimes a\big)=
{\rm det}(A){\rm det}\big(1+A^{-1}a\otimes a\big)={\rm det}(A)\Big(1+(A^{-1}a,a)\Big),
\end{eqnarray*}
since ${\rm tr}(A^{-1}a\otimes a)=(A^{-1}a,a)$ and  $c_k={\rm tr}\left(\textstyle \bigwedge ^{k}\big(A^{-1}a\otimes a\big)\right)=0$ for all $k>1$. To verify the last statement, by \eqref{tr(ext)^k} it is sufficient to verify that ${\rm tr}\big(D^k\big)= \big({\rm tr}\,D\big)^k$ for $D=A^{-1}a\otimes a$. Indeed, we have
${\rm tr}\,D=(A^{-1}a,a)$ and
\begin{equation}
 \label{tr(D^k)=tr(D)^k}
D^k=(A^{-1}a,a)^{k-1}D \quad \text{therefore},\quad {\rm tr}\big(D^k\big)= \big({\rm tr}\,D\big)^k.
\end{equation}
\qed\end{pf}
If we take $A=C(\lambda)$ we will get
\begin{equation}
\label{(C^{-1}(lam)a,a).2} 
1+(C(\lambda)^{-1}a,a)=\frac{{\rm det}\big(C(\lambda)+a\otimes a\big)}{{\rm det}\,C(\lambda)}.
%=\frac{{\rm det}(I+\gamma(y_1,y_2,\dots,y_m))}{{\rm det}(I+\gamma(y_2,\dots,y_m))}. 
\end{equation}

%by \eqref{Delta(f_k)k}  

%\commA{end )}
%
\begin{thm}
\label{t.m}
Let $C$ be defined by \eqref{C=Gram} and $a,\lambda\in \mathbb C^n$,
then
\begin{equation}
\label{m}
1+\big(
C(\lambda)^{-1}
a,a\big)=
\frac{{\rm det}\big(I_m+\gamma(y_1,y_2,\dots,y_m)\big)}{{\rm det}\big(I_{m-1}+\gamma(y_2,\dots,y_m)\big)}
=1+\Delta(y_1,y_2,\dots,y_m),
\end{equation}
where  $y_k$ for $1\leq k\leq m$ are  defined by \eqref{x_k,y_r=} and  $\Delta(y_1,y_2,\dots,y_m)$
is defined by \eqref{Delta(f_k)k}.
 \iffalse
% as follows
%\begin{equation}
%\label{y_k=} 
%y_k=y_k^{(n)}=\left(\frac{a_{kn}}{\sqrt{\lambda_k}}\right)_{k=1}^n.
\end{equation} 
\fi
\end{thm}
\begin{pf}
%
%\commA{14.10.23, small 1127-66}
%
{\small By Lemma~\ref{l.(A^{-1}a,a)}
it is sufficient to show that
\begin{eqnarray*}
&&
 {\rm det}\big(C(\lambda)\!+\!a\otimes a\big)\!=\!\Big(\prod_{k=1}^n\lambda_k\Big){\rm det}\big(I+\gamma(y_1,\dots,y_m)\big),\\
&& 
{\rm det}\,C(\lambda)\!=\!\Big(\prod_{k=1}^n\lambda_k\Big){\rm det}\big(I\!+\!\gamma(y_2,\dots,y_m)\big).
\end{eqnarray*}
Indeed, we have
\begin{eqnarray*}
&&
C(\lambda)+a\otimes a=\gamma(g_1,\dots,g_n)+{\rm diag}(\lambda_k)_{k=1}^n+(a_{1k}a_{1r})_{k,r=1}^n=\\
&&
\big((g_k,g_r)+a_{1k}a_{1r}\big)_{k,r=1}^n+{\rm diag}(\lambda_k)_{k=1}^n\stackrel{\eqref{X(3n)}}{=}\big((x_k,x_r)\sqrt{\lambda_k\lambda_r}\big)_{k,r=1}^n+\\
&&
{\rm diag}(\lambda_k)_{k=1}^n={\rm diag}(\sqrt{\lambda_k})_{k=1}^n\Big(I+\gamma(x_1,\dots,x_n)\Big){\rm diag}(\sqrt{\lambda_k})_{k=1}^n.\\
 &&
\text{Therefore,}\quad {\rm det}\,\big(C(\lambda)+a\otimes a\big)=
\Big(\prod_{k=1}^n\lambda_k\Big)
{\rm det}\big(I+\gamma(x_1,\dots,x_n)\big)\stackrel{\eqref{d(X^*X)=d(XX^*)}}{=} \\
&&
\Big(\prod_{k=1}^n\lambda_k\Big){\rm det}\big(I+\gamma(y_1,\dots,y_m)\big).
\end{eqnarray*}
Further, 
\begin{eqnarray*}
&&
{\rm det}\,C(\lambda)={\rm det}\Big(\gamma(g_1,\dots,g_n)+{\rm diag}(\lambda_k)_{k=1}^n\Big)
\stackrel{\eqref{x_k,y_r}}{=}
\\
&&
{\rm det}\,\Big({\rm diag}(\sqrt{\lambda_k})_{k=1}^n 
\big(I+\gamma(\bar{x}_1,\dots,\bar{x}_n)\big)
{\rm diag}(\sqrt{\lambda_k})_{k=1}^n\Big)=\\
&&
\Big(\prod_{k=1}^n\lambda_k\Big){\rm det}\,
\big(I+\gamma(\bar{x}_1,\dots,\bar{x}_n)\big)\stackrel{\eqref{d(X^*X)=d(XX^*)}}{=}
\Big(\prod_{k=1}^n\lambda_k\Big){\rm det}\,
\big(I+\gamma(y_2,\dots,y_m)\big).\quad \Box
\end{eqnarray*}
}
%\qed
\end{pf}

\section{Some estimates}
\label{s.estim}
We give  the well known estimates (see, e.g.,
%for example, 
\cite{BecBel61}, Chap. I, \S 52)
\begin{equation}
\label{A.min1}
 \min_{x\in{\mathbb
R}^n}\Big(\sum_{k=1}^na_kx_k^2\mid \sum_{k=1}^n x_k=1\Big)=
\Big(\sum_{k=1}^n \frac{1}{a_k}\Big)^{-1},\quad a_k>0,\,x_k\in \mathbb R.
\end{equation}
We will also use  the same estimate in a slightly different
form:
\begin{equation}
\label{A.min2}
 \min_{x\in{\mathbb
R}^n}\Big(\sum_{k=1}^na_kx_k^2\mid \sum_{k=1}^n x_kb_k=1\Big)=
\Big(\sum_{k=1}^n \frac{b_k^2}{a_k}\Big)^{-1}.
\end{equation}
The minimum is reached  for
$x_k=\frac{b_k}{a_k}\Big(\sum_{k=1}^n\frac{b_k^2}{a_k}\Big)^{-1}$.

%\label{s.1.4.1}
\begin{lem} [\cite{Kos04,Kos_B_09}]
\label{l.1.min}
For a 
 %strictly 
positive operator $A$, satisfying $(Af,f)>0,\\\,f\not=0$, acting on ${\mathbb R}^n$ and a vector
$b\in{\mathbb R}^n\backslash\{0\}$ we have
\begin{equation}
\label{A.min3} \min_{x\in{\mathbb R}^n}\{(Ax,x)\mid
(x,b)=1\}=\frac{1}{(A^{-1}b,b)}.
\end{equation}
The minimum is reached  for
$x=\frac{1}{(A^{-1}b,b)}A^{-1}b.$ 
\end{lem}
Lemma \ref{l.1.min} is a direct generalization of \eqref{A.min2}.
\iffalse
\begin{pf} 
We can study the corresponding 
Lagrange function
\begin{equation*}
 F(\lambda,t_1,\dots,t_n)\!=\!(At,t)-\lambda(b,t)
\!=\!\sum_{k=1}^na_{kk}t_k^2+2\sum_{1\leq k<r\leq n}a_{kr}t_kt_r-
\lambda\sum_{k=1}^nb_kt_k.
\end{equation*}	
Differentiating, we get
\begin{equation}
\label{dF-dt}
\frac{\partial F}{\partial t_k}
=2a_{kk}t_k+2\sum_{r=1,r\not=r}^na_{kr}t_r-
\lambda b_k=2\sum_{r=1}^na_{kr}t_r-
\lambda b_k=0.
\end{equation}
Relation \eqref{dF-dt}  means that $2At=\lambda b$, hence $t=\frac{\lambda}{2} A^{-1}b$.
Multiplying both part of \eqref{dF-dt} by $t_k$ and summing we get
\begin{equation}
\label{(At,t)=1}
(At,t)=\frac{\lambda}{2}(b,t)=\frac{\lambda}{2}.
\end{equation}
Substituting here  $t=\frac{\lambda}{2} A^{-1}b$ we receive
$
\frac{\lambda}{2}(A^{-1}b,b)=1
$
or, finally, $(At,t)=\frac{1}{(A^{-1}b,b)}$ by \eqref{(At,t)=1}.
\qed\end{pf}
%%
\fi
%%%%%%%%%%%%%%%%%%%%
\begin{lem} [\cite{Kos-m-Arx23}]
\label{l.A.min.H}
For a strictly positive operator $A$ on an infinite-dimensi\-onal Hilbert space $H$ and a vector $b\in H\backslash\{0\}$ 
such that $b\in D(A^{-1})$, {\rm where $D(B)$ is the domain  of the definition of an operator $B$},
we have
\begin{equation}
\label{A.min.H} 
\min_{x\in H}\{(Ax,x)\mid
(x,b)=1\}=\frac{1}{(A^{-1}b,b)}.
\end{equation}
The minimum is reached  for
$x_0=\frac{1}{(A^{-1}b,b)}A^{-1}b.$ 
\end{lem}
%%%
\begin{pf}
Consider a new scalar product in $H$ defined as follows:
\begin{equation}
\label{(.,.)new.H}
(f,g)_A:=(Af,g),\quad f,g\in H.
\end{equation} 
Since  $$(Ax,x)=(x,x)_A=\Vert x \Vert^2_A\quad \text{and}\quad 1=(b,x)=(A^{-1}b,x)_A,
$$ the minimum $\Vert x \Vert^2_A$ will be achived on the vector $x_0=sA^{-1}b$ generating hyperplane $1=(A^{-1}b,x)_A$ and lying on this hyperplane. We get
$$
1=(b,sA^{-1}b),\quad \text{therefore}\quad s=\frac{1}{(A^{-1}b,b)},\quad x_0=\frac{1}{(A^{-1}b,b)}A^{-1}b.
$$
Finally, we get
$
%\hskip 4 cm
(Ax_0,x_0)=\frac{1}{(A^{-1}b,b)}.
%\hskip 5 cm\qed
$
\qed\end{pf}
\begin{coex}
\label{cx.b-in-D(A)}
For a positive 
definite
operator $A={\rm diag}(\lambda_k)_{k=1}^\infty$ in $l_2(\mathbb N)$ where $\lambda_k=\frac{1}{k}$  and $b=(b_k)_{k\in \mathbb N}\in l_2(\mathbb N)$ with $b_k=\frac{1}{k}$ we have $b\not\in D(A^{-1})$, since $(A^{-1}b)_k\equiv1$ for all $k\in \mathbb N$ hence, $A^{-1}b\not\in l_2(\mathbb N)$. In this case $\frac{1}{(A^{-1}b,b)}=0$. Indeed, for the corresponding projections $A_n,\,\,b_n$ on $\mathbb R^n$ we have 
$$(A_n^{-1}b_n,b_n)=\sum_{k=1}^n\frac{1}{k}\to \infty.
$$
\end{coex}

%\newpage
\section{Application}
\subsection{The general idea}
In the concrete examples considered in 
\cite{Kos92}--\cite{KosMor-Arx23}
the possibility to approximate a lot of functions in $L^\infty (X_m,\mu)$ using Lemma \ref{l.1.min},   follows from the fact that
\begin{equation}
\label{to-inft}
\lim_{n\to\infty}(C_n(\lambda)^{-1}a_n,a_n) =\infty.
\end{equation}
By Theorem~\ref{t.m} we have
\begin{equation}
\label{main-form} 
\big(C_n(\lambda)^{-1}a_n,a_n\big)\!=\!\Delta(y_1^{(n)},y_2^{(n)},\dots,y_m^{(n)})\!=\!\frac{{\rm det}\big(I_m+\gamma(y_1^{(n)},y_2^{(n)},\dots,y_m^{(n)})\big)}{{\rm det}\big(I_{m-1}+\gamma(y_2^{(n)},\dots,y_m^{(n)})\big)}-1.
\end{equation}
Finally, by Lemma 
\ref{l.min=proj.m} and 
Lemma~\ref{l.det/det.m}, proved in \cite{Kos-m-Arx23}, we have
\begin{equation*}
\lim_{n\to\infty}\!
\frac{{\rm det}\big(I_m+\gamma(y_1^{(n)},y_2^{(n)},\dots,y_m^{(n)})\big)}{{\rm det}\big(I_{m-1}+\gamma(y_2^{(n)},\dots,y_m^{(n)})\big)}
\!=\!\infty. 
\end{equation*}
%For convenience of the reader we formulate  
%
\begin{lem}[\cite{Kos-m-Arx23}]
\label{l.min=proj.m}
Let $f_r=(f_{rk})_{k\in \mathbb N},$
%$m\in \mathbb N$ and $(f_r)_{r=0}^{m}$ 
be $m+1$ infinite real vectors 
$0\leq r\leq m$
such that for all $\big(C_0,\dots,C_{m}\big)\in {\mathbb R}^{m+1}\setminus\{0\}$  holds
\begin{equation}
\label{norm=infty.m}
\sum_{r=0}^{m}C_rf_r\not\in l_2(\mathbb N),\quad\text{i.e.,} \quad
\sum_{k\in \mathbb N}\Big| \sum_{r=0}^{m}C_rf_{rk}\Big|^2=\infty.
\end{equation}
Denote by $f_r^{(n)}=(f_{rk})_{k=1}^n\in \mathbb R^n$ the {\rm projections} of the vectors $f_r$ on the subspace $\mathbb R^n$. Then for all 
$s$ with $0\leq s\leq m$
\begin{equation}
\label{final.m}
\frac{\Gamma(f_0,f_1,\dots,f_{m})}
{\Gamma(f_0,\dots,\hat{f_s},\dots,f_{m})}:=
\lim_{n\to\infty}\frac{\Gamma(f_0^{(n)},f_1^{(n)}\dots,f_{m}^{(n)})}{\Gamma(f_0^{(n)},\dots,\widehat{f_s^{(n)}},\dots,f_{m}^{(n)})}=\infty,
\end{equation}
where $\hat{f_s}$ means that the vector $f_s$ is absent and $\Gamma(f_0,f_1,\dots,f_{m})$ is {\rm the Gram determinant}.
\end{lem}

\begin{lem}[\cite{Kos-m-Arx23}]
\label{l.det/det.m} 
Let we have $m+1$ real vectors $(f_k)_{k=0}^m$  such that
$\sum_{k=0}^mC_kf_k\\\not\in l_2(\mathbb N)$ for {\rm any  nontrivial combination} $(C_k)_{k=0}^m$.
Then for any $s,\,0\leq s\leq m$ %holds
\begin{equation}
\label{det/det.m}
\frac{{\rm det}
\big(I_{m+1}+\gamma(f_0,\dots, f_m)\big)}
{{\rm det}\big(I_{m}\!+\!\gamma(f_0,\dots,\hat{f_s},\dots,f_m)\big)}\!=\!
\lim_{n\to\infty}\!
\frac{
{\rm det}\big(I_{m+1}+\gamma(f_0^{(n)},\dots, f_m^{(n)})\big)}
{{\rm det}\big(I_{m}\!+\!\gamma(f_0^{(n)},\dots,\widehat{f_s^{(n)}},\dots,f_m^{(n)})\big)}
\!=\!\infty.
\end{equation}
Here $I_m\!=\!{\rm diag}(1,\dots,1)\in{\rm Mat}(m,\mathbb R,)$
and $\gamma(f_0,\dots, f_m)$ is  {\rm the Gram matrix}.
%for vectors $(f_k)_{k=0}^m$.
\end{lem}
\begin{pf}
The proof follows from Lemma~\ref{l.min=proj.m} and  \eqref{detC-LI}.
\qed\end{pf}
\begin{rem}
\label{r.h_n}
We note that $\frac{\Gamma(f_0,f_1,\dots,f_m)}
{\Gamma(f_1,\dots,f_m)}$ is the square of the {\it height} of the  
{\it parallelotope}
generated by the vectors $f_0,f_1,\dots,f_m\in \mathbb R^{m+1}$, see Lemma~\ref{l.d(f,v_n)}.
\end{rem}
\index{parallelotope}
\index{height}

%\subsection{One example of the group $B_0^{\mathbb N}$}
%\newpage
\subsection{The Ismagilov conjecture}

\label{Ism.comj}
To construct the regular representation for an infinite-dimensi\-onal
group $G$, first we should find some larger topological group $\widetilde{G}$
and a measure $\mu$ on $\widetilde{G}$ such that $G$ is a dense
subgroup in $\widetilde{G},$
and the measure is  right or left $G$-quasi-invariant, i.e.,
$\mu^{R_t}\sim\mu$ for all $t\in
G,$ (or $\mu^{L_s}\sim\mu$ for all $s\in G$), here $\sim$ means {\it equivalence}, for details see \cite{Kos_B_09}. 
We use notation $\mu^f(\Delta)=\mu\big(f^{-1}(\Delta)\big)$ for $f:X\to X$, where $\Delta$ is some measurable set in $X$.
Consider the right and the left actions $R_t,L_s$ of the group
$G$ on $\widetilde
G$ defined below:
$$
R_tx=xt^{-1},\quad L_sx=sx,\quad t,s\in G,\,\,x \in \widetilde
G.
$$
Denote by $\mu^{R_t},\,\,\mu^{L_s}$ the images of the measure
$\mu$ under the map $R_t,L_s:\widetilde
G\to \widetilde G$.
The right and left
representations $T^{R,\mu},T^{L,\mu}:G\rightarrow U(L^2(\widetilde
G,\mu))$ are naturally defined in the Hilbert space $L^2(\widetilde
G,\mu)$ by the following formulas:
\begin{equation}
\label{T(R,b)}
(T^{R,\mu}_tf)(x)=(d\mu(xt)/d\mu(x))^{1/2}f(xt),
\end{equation}
\begin{equation}
\label{T(L,b)}
(T^{L,\mu}_sf)(x)=(d\mu(s^{-1}x)/d\mu(x))^{1/2}f(s^{-1}x).
\end{equation}
\index{measure!quasi-invariant}
The right regular representation of infinite-dimensional
groups can be irreducible if no left actions are {\it admissible} for
the measure $\mu$, i.e., when $\mu^{L_t}\perp \mu$ for all $t\in
G\backslash{\{e\}}$. In this case a von Neumann algebra
${\mathfrak A}^{T^{L,\mu}}$ generated by the left regular
representation $T^{L,\mu}$  is trivial.  More precisely:

\begin{co}
[Ismagilov, 1985] \label{co.Ism-art} The right regular
representation \vskip -0.2cm
$$T^{R,\mu}:G\rightarrow U(L^2(\widetilde G,\mu))$$
%\noindent
is irreducible if and only if \par 1) $\mu^{L_t}\perp
\mu\quad\text{for all}\quad t\in G\backslash{\{e\}},\,\,$ (where $\perp$ stands
for singular),
\par
2) the measure $\mu$ is $G$-ergodic.
\end{co}
\index{measure!ergodic}
Conditions 1) and 2) are the necessary conditions of the irreducibility. The problem is to prove that they are sufficient ones too.
\begin{rem}
\label{r.Ism.conj}
 This conjecture was expressed by Rais Salmanovich Ismagilov in his referee report of the author's PhD Thesis, 1985. It was verified for a lot of particular cases. In the general case, it is an open problem.
In the case of a finite field ${\mathbb F}_p$ we need some additional conditions for the irreducibility  
\cite{Kos16F(p).ar}.
\end{rem}
\iffalse
For symmetric groups and inductive limits of classical compact groups there are
analogs of regular representations of another type (with a well-developed harmonic analysis)
see \cite{BorOlsh05,KerOlshVerh04}.}
\fi

%
\index{conjecture!Ismagilov}
\subsection{Group  $B_0^{\mathbb N}$, arbitrary  mesure $\mu$}
\label{sec.reg.B^N.mu}

Let $B_0^{\mathbb N}$ be the group of finite real
upper-triangular matrices with unities on the principal diagonal
and let $B^{\mathbb N}$ be the group of all such matrices (not
necessarily finite):
$$
\aligned B_0^{\mathbb N} & =\{I+x=I+\sum_{k<n}x_{kn}E_{kn} \mid x
\,\,{\text {is\,\,finite}}\},       \\
B^{\mathbb N} & =\{I+x=I+\sum_{k<n}x_{kn}E_{kn} \mid x\,\,\text
{is\,\, arbitrary}\}.
\endaligned
$$
Let $\mu$ be an arbitrary probability measure on the group
$B^{\mathbb N}$. If $\mu^{R_{t}}\sim\mu $ and
$\mu^{L_{t}}\sim\mu\quad\text{for all}\quad t\in B_0^{\mathbb N}$,  an
analogue of the right $T^{R,\mu}$ and the left $T^{L,\mu}$ regular
representations of the group $B_0^{\mathbb N}$, i.e.,
$T^{R,\mu},\,\,T^{L,\mu}:B_0^{\mathbb N}\rightarrow U(H_\mu)$ are
defined in the space $H_\mu=L^2( B^{\mathbb N},\mu)$ by \eqref{T(R,b)} and  \eqref{T(R,b)}.
\iffalse
In this case the Ismagilov conjecture has the following form
\begin{co}
[Ismagilov, 1985] \label{co.Ism.B^N} The right regular
representation \vskip -0.2cm
$$T^{R,\mu}:B_0^{\mathbb N}\rightarrow U(L^2(B^{\mathbb N},\mu))$$
%\noindent
of the group $B_0^{\mathbb N}$ is irreducible if and only if 
\end{co}
\begin{equation}
\label{T(R,b)}
(T^{R,\mu}_t f)(x)  =(d\mu(xt)/d\mu(x))^{1/2}f(xt),
\end{equation}
\begin{equation}
 \label{T(L,b)}
(T^{L,\mu}_t f)(x) 
=(d\mu(t^{-1}x)/d\mu(x))^{1/2}f(t^{-1}x).
\end{equation}
%
%%%
\fi

For the generators $A^{R,\mu}_{kn}$ ($A^{L,\mu}_{kn}$) of the
one-parameter groups $I+tE_{kn},\,t\in {\mathbb R},\,k<n,$
corresponding to the right $T^{R,\mu}$ (respectively the left
$T^{L,\mu}$) regular representation we have the following
formulas:

\begin{eqnarray}
\label{56.4}  
&&
A^{R,\mu}_{kn}  =
\frac{d}{dt}T^{R,\mu}_{I+tE_{kn}}\vert_{t=0}=
\sum_{r=1}^{k-1} x_{rk}D_{rn}(\mu)+D_{kn}(\mu),
\\
\label{56.5} 
&&
A^{L,\mu}_{kn}=
%\left.
\frac{d}{dt}T^{L,\mu}_{I+tE_{kn}}
%\right
\vert_{t=0}=
-(D_{kn}(\mu)+\sum_{m=n+1}^\infty x_{nm}D_{km}(\mu)),
\end{eqnarray}
where $D_{kn}(\mu)=\displaystyle \frac{\partial}{\partial
x_{kn}}+
%\left. 
\frac{d}{dt} \bigg(
\frac{d\mu(x(I+tE_{kn}))}{d\mu(x)}\bigg) ^{1/2}%\right
\vert_{t=0}.$
For an arbitrary  product measure  $\mu =\otimes_{k<n}\mu_{kn}$,
we have
\begin{equation}
\label{D(kn)} 
D_{kn}(\mu)=\frac{\partial}{\partial x_{kn}}+
\frac{\partial}{\partial x_{kn}}\Big(
\ln\mu_{kn}^{1/2}(x_{kn})\Big) ,
\end{equation}
where we write $d\mu_{kn}(x)=\mu_{kn}(x)dx,\,x\in {\mathbb R}.$

\subsubsection{Group $B_0^{\mathbb N}$, Gaussian centered mesure}
\label{sec.reg-B(N),(b)}
See details in \cite[Ch. 2.1]{Kos_B_09}. 
Let us define the Gaussian product-measure $\mu_b$ on the group $B^{\mathbb
N}$ in the following way:
\begin{equation}
\label{mu(b)} 
d\mu_b(x) =
\otimes_{k<n}(b_{kn}/\pi)^{1/2}\exp(-b_{kn}x_{kn}^2)dx_{kn}
=\otimes_{k<n}d\mu_{b_{kn}}(x_{kn}),
\end{equation}
where $b=(b_{kn})_{k<n}$ is some set of positive numbers. In
this case we have 
%(see, for example, \cite{KosZek00}, formulas (6)
%and (7), or \cite[Ch. 2.2.1, formula (2.8)]{Kos_B_09}).
\begin{equation}
\label{A_{kn}}
A^{R,\mu}_{kn}  =
%\left.
\frac{d}{dt}T^{R,\mu}_{I+tE_{kn}}
%\right
\vert_{t=0}=
\sum_{r=1}^{k-1} x_{rk}D_{rn}+D_{kn},\quad
D_{kn}=\frac{\partial}{\partial x_{kn}}- b_{kn}x_{kn},
\end{equation}

It turns out that the measure $\mu_b$ is always $B_0^{\mathbb
N}$-right-quasi-invariant.
%(Lemma \ref{l.mu-R-eqv}). 
Therefore, we can
construct a family of analogues of the right $T^{R,\mu_b}$ and the left
$T^{L,\mu_b}$ (if the measure $\mu_b$ is $B_0^{\mathbb N}$-left-quasi-invariant)
regular representations of the group $B_0^{\mathbb N}$ in the space
\index{representation!regular}
$
%H(b)=\!
L_{\,2}(B^{\mathbb N},\mu_b).
$
They are defined by \eqref{T(R,b)} and  \eqref{T(L,b)}.
\begin{thm}
[\cite{Kos92,Kos_B_09}]
\label{T^(R,b)-irr}
The right regular representation $T^{R,\mu_b}$ of the group
$B_0^{\mathbb N}$ is irreducible if and only if 
%$\mu_b^{L_t}\perp
%\mu_b\quad\text{for all}\quad t\in B_0^{\mathbb N}%\backslash{\{e\}}.$
\par 1) $\mu^{L_t}\perp
\mu\quad\text{for all}\quad t\in B_0^{\mathbb N}\backslash{\{e\}},$
\par
2) the measure $\mu$ is $B_0^{\mathbb N}$-ergodic.
\end{thm}
\begin{df}
\label{d.erg}
Let $\alpha:G\rightarrow {\rm
Aut}(X)$  be a {\it measurable action} of a group $G$  on a
measurable space $(X,%{\mathfrak B},
\mu)$.
Recall that the probability measure $\mu$ on some $G$-space $X$ is called {\it ergodic} if any function
$f\in L^1(X,\mu)$ with property $f(\alpha_t(x))=f(x)$ a.e. (almost everywhere) ${\rm mod}\,\mu$ is constant.
\end{df}
\begin{lem}
[\cite{Kos_B_09}, Lemma 2.1.6]
\label{l.orth-B_0^N}
%\label{l.56.3}
We have $\mu_b^{L_t}\perp\mu_b\quad\text{for all}\quad  t\in B_0^{\mathbb
N}\backslash{e}$
if and only if
\begin{equation}
\label{S^L_{kn}(b)} 
%\Leftrightarrow\,\,
S_{kn}^L(\mu_b)=\sum_{m=k+1}^\infty\frac{b_{km}}{b_{nm}}=\infty
\quad\text{for all}\quad k\!<\!n.
\end{equation}
\end{lem}
{\it Idea of the proof of irreducibility, for details see} \cite{Kos92,Kos_B_09}.
The conditions 1) and 2 are necessary conditions of the irreducibility of the representation $T^{R,\mu_b}$.
We show that they are also a sufficient ones. Let $\mathfrak A
%^{R,\mu_b}
( B_0^{\mathbb N})$ be a von Neumann algebra generated by the representation  $T^{R,\mu_b}$:
\begin{equation}
\label{A-von-Neum.1}
\mathfrak A
(B_0^{\mathbb N})=\Big(T^{R,\mu_b}_t\mid t\in B_0^{\mathbb N}\Big)''.
\end{equation}
To prove the irreducibility, it is sufficient to show that 
$U_{kn}(t)
\in \mathfrak A(B_0^{\mathbb N})$
for all $k,n\in \mathbb N,\,k<n$, where $U_{kn}(t)=e^{itx_{kn}}$. In this case we have 
\begin{equation}
 \label{L(inf)-in-A}
L^\infty(B^{\mathbb N},\mu_b)\!\subset\! \mathfrak A
(B_0^{\mathbb N})\quad\text{hence,}\quad
%\\&&
 \big(\mathfrak A
(B_0^{\mathbb N})\big)'\!\subset\!
\Big(L^\infty(B^{\mathbb N},\mu_b)\Big)'=L^\infty(B^{\mathbb N},\mu_b),
\end{equation}
since the algebra $L^\infty(B^{\mathbb N},\mu_b)$ is {\it maximal abelian}.
Let now some bounded operator $A$ commute with the representation $[T^{R,\mu_b}_t,A]=0$ for all $t\in B_0^{\mathbb N}$. Then by \eqref{L(inf)-in-A},
%an operator 
$A\in L^\infty(B^{\mathbb N},\mu_b)$ , i.e, $A$ is a multiplication operator on some function $a\in L^\infty(B^{\mathbb N},\mu_b)$. The commutation $[T^{R,\mu_b}_t,a]=0$ implies
$a(xt)=a(x)$ a.e. {\rm mod}\,$\mu_b$. By ergodicity of the measure $\mu_b$ on $B^{\mathbb N}$ we conclude that $a(x)=const$ hence $A=CI$, i.e, representation $T^{R,\mu_b}$ is irreducible.

To illustrate the approximation we show here only that $e^{itx_{12}}\in \mathfrak A
(B_0^{\mathbb N})$, or $x_{12}\,\,\eta\,\, \mathfrak A
(B_0^{\mathbb N})$, i.e., that operator $x_{12}$ is  {\it affiliated} with an algebra $\mathfrak A
(B_0^{\mathbb N})$. 
\begin{df}
\label{A-eta-M} Recall 
%(see \cite{Dix69W}) 
that, a  not
necessarily bounded self-adjoint operator $A$ in a Hilbert space
$H$, is said to be {\it affiliated} with a von Neumann algebra $M$
of operators in this Hilbert space $H$ if $e^{itA}\in M$ for all
$t\in{\mathbb R}$. One writes $A\,\,\eta\,\,M$.
\end{df}
We show that the operator $x_{12}$ can be approximated in the {\it strong resovent sense}
by the linear combinations of the following operators $A_{1k}A_{2k},\,\,k\geq 3$. By \eqref{A_{kn}} we get
\begin{equation}
 A_{1k}A_{2k}=D_{1n}(x_{12}D_{1k}+D_{2k})=x_{12}D_{1k}^2+D_{1k}D_{2k}\,\,k\geq 3.
\end{equation}
By \cite{Kos_B_09}, Lemma 2.1.9, the convergence 
$
\sum_{k=N_1}^{N_2}t_k  A_{1k}A_{2k}\to x_{12}
$
holds if and only if $S_{12}^L(\mu_b)=\sum_{k=3}^\infty\frac{b_{1k}}{b_{2k}}=\infty$. And this is precisely the condition of orthogonality $\mu_b^{L_t}\perp\mu_b$, see Lemma~\ref{l.orth-B_0^N}.

We give here more conceptual proof of this fact. Using the appropriate Fourier transform $F_2$ in the variables $(x_{1k},x_{2k})_{k=3}^\infty$ see details in \cite[Section 2.1.3, formula (2.15)]{Kos_B_09}  we get $F_2(D_{1k})=y_{1k},\,\,F_2(D_{2k})=y_{2k},\,\,k\geq 3$ therefore,
\begin{equation}
F_2(A_{1k}A_{2k})=x_{12}y_{1k}^2+ y_{1k}y_{2k}.
\end{equation}
The corresponding measure $\mu_{1/4b}(y)$ in  variables $(y_{1k},y_{2k})_{k=3}^\infty$ 
is defined by 
\begin{equation}
\label{mu-tilde}
d\mu_{1/4b}(y) =
\otimes_{k=1}^2\otimes_{n=3}^\infty
\sqrt{\frac{1}{4b_{kn}\pi}}
\exp\Big(-
\frac{y_{kn}^2}{4b_{kn}}\Big)dy_{kn}
=\otimes_{k=1}^2\otimes_{n=3}^\infty d\mu_{1/4b_{kn}}(y_{kn}).
\end{equation}
The corresponding canonical measure $\mu_{1/2}(z)$ is as follows:
\begin{equation}
\label{mu-stand.2}
d\mu_{1/2}(z) =
\otimes_{k=1}^2\otimes_{n=3}^\infty
\sqrt{\frac{1}{2\pi}}
\exp\Big(-
\frac{z_{kn}^2}{2}\Big)dz_{kn}
=\otimes_{k=1}^2\otimes_{n=3}^\infty d\mu_{1/2}(z_{kn}).
\end{equation}
In the {\it canonical coordinates} $z_{kn}$ the expression $F_2(A_{1k}A_{2k})$ will have the following form
$$
x_{12}2b_{1k}z_{1k}^2+ 2\sqrt{b_{1k}b_{2k}}z_{1k}z_{2k}=2\sqrt{b_{1k}}\big(x_{12}z_{1k}^2+ a_kz_{1k}z_{2k}\big),\quad a_k=\sqrt{b_{2k}/b_{1k}}.
%\sqrt{\frac{b_{2k}}{b_{1k}}}.
$$
Let us denote by $\langle f_n\mid n\in {\mathbb
N}\rangle$ the {\it closure of the linear space}
generated by the set of
vectors $(f_n)_{n\in{\mathbb N}}$ in a Hilbert space $H.$ 
\begin{lem}
\label{l.x_12-can-coord} 
Set $f_0=x_{12}, f_k=x_{12}z_{1k}^2+ a_kz_{1k}z_{2k}$. 
We have 
%
%\commA{11.10.23, expalin notations}
%
%{\color{red}
$f_0\in 
\langle  f_k\mid k\geq 3 \rangle$
if and only if $\sum_{k=3}^\infty\frac{1}{a_k^2}=\infty$.
\end{lem}
\begin{pf}
Consider the hyperplane $V_n$ generated by $n$  vectors $f_3,\dots, f_{n+3}$. 
By Lemma~\ref{l.d(f,v_n)} we have
\begin{equation}
\label{d(f,v_n).1}
d^2(f_0,V_n)=
%\frac{G(f_0,f_1,f_2,\dots, f_n,)}{G(f_1,f_2,\dots, f_n)}.
\frac{\Gamma(f_0,f_3,f_4,\dots, f_{n+3})}{\Gamma(f_3,f_4,\dots, f_{n+3})}.
\end{equation}
Further,
\begin{equation}
\label{gam(0,3,..,n+3)}
\gamma(f_0,f_3,f_4,\dots,f_{n+3})
= \left(
\begin{array}{cccc}
	1&1          &...&1\\
	1          &1+a^2_3&...&1\\
	&&...&\\
	1          &          1&...&1+a^2_{n+3}
	\end{array}
	\right)
	\end{equation}
and
\begin{equation}
\label{gam(3,..,n+3)}
\gamma(f_3,f_4,\dots,f_{n+3})
= \left(
\begin{array}{cccc}
	1+a^2_3&1          &...&1\\
	1          &1+a^2_4&...&1\\
	&&...&\\
	1          &          1&...&1+a^2_{n+3}
	\end{array}
	\right).
	\end{equation}
Finally, by \eqref{det(C).1} we get	
\begin{eqnarray*}
&&
d^2(f_0,V_n)=\frac{{\rm det}\big(\gamma(f_0,f_3,f_4,\dots,f_{n+3})\big)}{{\rm det}\big(\gamma(f_3,f_4,\dots,f_{n+3})\big)}\stackrel{\eqref{det(C).1}}{=}
\frac{\Big(\prod_{k=3}^{n+3}a^2_k\Big)}{\Big(\prod_{k=3}^{n+3}a^2_k\Big)
\Big(1+\sum_{k=3}^{n+3}\frac{1}{a^2_k}\Big)}
\\
&&
\hskip 1.8cm
=
\Big(1+\sum_{k=3}^{n+3}\frac{1}{a^2_k}\Big)^{-1}.\hskip 7.5
cm \Box
\end{eqnarray*}
\end{pf}

\subsection{Koopman's representation}
Let $\alpha:G\rightarrow {\rm Aut}(X)$ be a measurable action of a
group $G$  on a measurable space $(X,\mu)$ with
$G$-quasi-invariant measure $\mu$, i.e, $\mu^{\alpha_t}\sim\mu$
for all $t\in G$. With these date one can associate the
representation $ \pi^{\alpha ,\mu,X}:G\rightarrow U(L^2(X,\mu)),
$ by the following formula:
\index{group!action!measurable}
\begin{equation}
\label{Rep(G,X)-pi}
(\pi^{\alpha,\mu,X}_tf)(x)=(d\mu(\alpha_{t^{-1}}(x))/d\mu(x))^{1/2}f(\alpha_{t^{-1}}(x)),\quad
f\in L^2(X,\mu).
\end{equation}
In the case of an invariant measure this representation called {\it Koopman's representation}.
%, see  \cite{Koo31}. 
We keep the same name for representation \eqref{Rep(G,X)-pi}.
\index{represenation!Koopman's}
\index{subgroup!centarlizer}
\index{represenation!Koopman's}
The following conjecture is a natural generalization of Ismagilov's conjecture.
\begin{co} 
%[Kosyak, \cite{Kos94,Kos02.3}]
\label{co.G-Ism-5}
The representation \eqref{Rep(G,X)-pi}
%$$\pi^{\alpha ,\mu,X}:G\rightarrow U(L^2(X,\mu))$$
is irreducible if and only if
\par 1) $\mu^g\perp \mu\quad\text{for all}\quad
g\in Z_{{\rm Aut}(X)}(\alpha(G))\backslash\{e\},\,\,$
\par 2) the measure $\mu$ is $G$-ergodic.
\end{co}
Here $Z_{G}(H)$ is a {\it centralizer} of the subgroup $H$ in the group $G$:
$
Z_{G}(H)=\{g\in G\mid \{g,a\}=e\,\,\forall a\in H\},
$
where $\{g,a\}=gag^{-1}a^{-1}$. 
In general, Conjecture~\ref{co.G-Ism-5} is false, our aim is to find when it holds, see the following section.

\subsection{Group ${\rm GL}_0(2\infty,{\mathbb R})$ acting on $m$ infinite rows}
Let us denote by ${\rm Mat}(2\infty,{\mathbb R})$ the space of all
real  matrices that are infinite in both directions:
\begin{equation}
\label{5.Mat(2inf,R)} {\rm Mat}(2\infty,{\mathbb R})=
\Big\{x=\sum_{ k,n\in{\mathbb Z}}x_{kn}E_{kn},\,\,x_{kn}\in
{\mathbb R}\Big\}.
\end{equation}

The group $G={\rm GL}_0(2\infty,{\mathbb R})=\varinjlim_{n,i^s}{\rm
GL}(2n+1,{\mathbb R})$ is defined as the inductive limit of the
general linear groups $G_n={\rm GL}(2n+1,{\mathbb R})$ with
respect to the {\it symmetric embedding} $i^s$:
\begin{equation}
\label{N.i^s} 
G_n\ni x\mapsto
i^s_{n+1}(x)=x+E_{-(n+1),-(n+1)}+E_{n+1,n+1}\in %
G_{n+1}.
\end{equation}
For a fixed natural number $m$, consider a $G$-space $X_m$  as the
following subspace of the space ${\rm Mat}(2\infty,{\mathbb R})$:
\begin{equation}
\label{5.X_m} X_m=\Big\{x\in {\rm Mat}(2\infty,{\mathbb
R})\mid x=\sum_{k=1}^m\sum_{n\in{\mathbb
Z}}x_{kn}E_{kn}\Big\}.
\end{equation}
The  group ${\rm GL}_0(2\infty,{\mathbb R})$ acts from the right
on the space $X_m.$ Namely, the right action of the group ${\rm
GL}_0(2\infty,{\mathbb R})$ is correctly defined on the space
$X_m$ by the formula $R_t(x)=xt^{-1},\,\,t\in G,\,\,x\in X_m$. We
define a Gaussian non-centered product measure
$\mu:=\mu^m:=\mu_{(b,a)}^m$
on the space $X_m:$
\begin{equation}
\label{5.mu^m}
\mu_{(b,a)}^m(x)=\otimes_{k=1}^m\otimes_{n\in{\mathbb
Z}}\mu_{(b_{kn},a_{kn})}(x_{kn}),
\end{equation}
where
\begin{equation}
\label{mu(b,a)} 
%\text{where}\quad 
d\mu_{(b_{kn},a_{kn})}(x_{kn})=\sqrt{\frac{ b_{kn}}{\pi}}e^{-b_{kn}(x_{kn}-a_{kn})^2}
dx_{kn}
\end{equation}
and $b=(b_{kn})_{k,n},\,\,b_{kn}>0,\,a=(a_{kn})_{k,n},\,a_{kn}\in
{\mathbb R},\,1\leq k\leq m,\,n\in {\mathbb Z}.$ Define the unitary 
representation $T^{R,\mu,m}$ of the group ${\rm
GL}_0(2\infty,{\mathbb R})$ on the space $L^2(X_m,\mu^m_{(b,a)})$
by the formula:
\begin{equation}
\label{T^(R,mu,m)}
(T^{R,\mu,m}_tf)(x)=\big(d\mu_{(b,a)}^m(xt)/d\mu_{(b,a)}^m(x)\big)^{1/2}f(xt),\,\,f\in L^2(X_m,\mu^m_{(b,a)}).
\end{equation}
Obviously, the {\it centralizer} $Z_{{\rm Aut}(X_m)}(R(G))\subset
{\rm Aut}(X_m)$ contains  the group $L({\rm GL}(m,{\mathbb R}))$,
i.e., the image of the group ${\rm GL}(m,{\mathbb R})$ with
respect to the left action $L:{\rm GL}(m,{\mathbb R})\rightarrow
{\rm Aut}(X_m),\,L_s(x)\!=\!sx,\,s\in {\rm GL}(m,{\mathbb
R}),\,x\in X_m.$  We prove the following theorem.
%for $m=3$. 
%
\begin{thm}
\label{5.t.irr} The representation $T^{R,\mu,m}\!:\!{\rm
GL}_0(2\infty,{\mathbb R})\!\rightarrow\! U\Big(L^2(X_m,\mu^m_{(b,a)})\Big)$ is
irreducible
%, for $m=3$, 
if and only if  
\begin{eqnarray*}
(i)&&(\mu^m_{(b,a)})^{L_{s}}\perp
\mu^m_{(b,a)}\quad\text {for all}\quad s\in {\rm GL}(m,{\mathbb
R})\backslash\{e\};\\
(ii)&&\text{the measure} \quad\mu^m_{(b,a)}\quad \text{is $G$-{\rm ergodic}}.
\end{eqnarray*}
\end{thm}
In {\rm \cite{Kos_B_09,KosJFA17}}
this result was proved for  $m\leq 2$. In  \cite{KosMor-Arx23} it was proved for $m=3$.
Note that conditions (i) and (ii) are necessary conditions for irreducibility.
\begin{rem}
Any Gaussian product-measure $\mu_{(b,a)}^m$ on $X_m$ is ${\rm
GL}_0(2\infty,{\mathbb R})$-right-ergodic \cite[\S 3, Corollary 1]{ShFDT67}, see Definition~\ref{d.erg}. For non-product-measures
this is not true in general.
\end{rem}
%%%

\subsubsection{Case $m=3$}
%
%For the case $m=2$, see details in \cite{Kos_B_09,KosJFA17}.
\begin{rem}
\label{r.irr-idea-G} (The idea of the proof of
irreducibility, see details in \cite{KosMor-Arx23}). Let us denote by ${\mathfrak A}^m$ the  {\it von
Neumann algebra} generated by the representation $T^{R,\mu,m}$, i.e., $
{\mathfrak A^m}=(T^{R,\mu,m}_t\mid t\in G)''. $ For
$\alpha\!=\!(\alpha_k)\!\in\!\{0,1\}^m$ define the von Neumann algebra
$L^\infty_\alpha(X_m,\mu^m)$ as follows:
$$
L^\infty_\alpha(X_m,\mu^m)\!=\!\Big(\exp(itB^\alpha_{kn})\mid
1\leq k\leq m,\,\,t\in {\mathbb R},\,\,n\in {\mathbb
Z}\Big)'',\,\,
$$
where $
B^\alpha_{kn}\!=\!\left\{\begin{array}{rcc}
 x_{kn},&\text{if}&\alpha_k=0\\
 i^{-1}D_{kn},&\text{if}&\alpha_k=1
\end{array}\right.$ and $D_{kn}=\partial/\partial x_{kn}-b_{kn}(x_{kn}-a_{kn})$.
\vskip 0.3 cm
{\bf The proof of the  irreducibility is based on  four
facts}:
\par 1) we can approximate by generators
$A_{kn}=A_{kn}^{R,m}=\frac{d}{dt}T^{R,\mu,m}_{I+tE_{kn}}\vert_{t=0}$
the set of operators $(B^\alpha_{kn})_{k=1}^m,\,n\!\in\!{\mathbb
Z}$ {\it for some} $\alpha\!\in\!\{0,1\}^m$ depending on the measure
$\mu^m$  using the orthogonality condition
$(\mu^m)^{L_{s}}\perp \mu^m$ for all $s\in {\rm GL}(m,{\mathbb
R})\backslash{\{e\}}$,
\par 2)
it is sufficient to  verify the approximation  only on the {\it cyclic vector} 
${\bf 1}(x)\!\equiv\! 1$, since the representation $T^{R,\mu,m}$ is {\it cyclic}, 
\par 3) the subalgebra $L^\infty_\alpha(X_m,\mu^m)$
is a {\it maximal abelian subalgebra} in ${\mathfrak A}^m$,
\par 4)  the measure $\mu^m$ is $G$-ergodic. 

Here the 
{\it generators}
$A_{kn}$ are given by the formulas:
$$
A_{kn}\!=\!\sum_{r=1}^{m}x_{rk}D_{rn},\quad k,n \in
{\mathbb Z},\quad\text{where}\quad D_{kn}=\partial/\partial x_{kn}-b_{kn}(x_{kn}-a_{kn}).
$$
\end{rem}
\begin{rem}
\label{r.shem} 
{\it Scheme of the proof.} We prove the irreducibility as follows
\begin{eqnarray}
\label{shem}
&&\left(\mu^{L_s}\perp \mu\,\,\, \text{for all} \,\,\, 
s\in {\rm GL}(3,\mathbb R)\setminus \{e\} \right)\Leftrightarrow
\left(\begin{smallmatrix}
\text{criteria}\\
\text{of}\\
\text{orthogonality}
\end{smallmatrix}\right)
\& \\
\nonumber
&&
\left(\begin{smallmatrix}
\text{Lemma~\ref{l.min=proj.m}
%{l.min=proj.3}
}\\
\text{about}\\
\text{three vectors}\, f,g,h\not\in l_2
\end{smallmatrix}\right)
\Rightarrow
\left(\begin{smallmatrix}
\text{some of}        &\Delta^{(1)},&\Delta_1\\
\text{the expressions}&\Delta^{(2)},&\Delta_2\\
\text{are divergent:} &\Delta^{(3)},&\Delta_3
\end{smallmatrix}
\right)
\Rightarrow
\text{irreducibility}, \\
\label{Si,si}
&&
\text{where}\quad \Delta^{(i)}:=\Delta(Y_i^{(i)},Y_j^{(i)},Y_k^{(i)}),\quad 
\Delta_i:=\Delta(Y_i,Y_j,Y_k), 
\end{eqnarray}
$\Delta(f,g,h)$ is defined by \eqref{Delta(f_k)k}, 
and $\{i,j,k\}$ is a cyclic permutation 
of $\{1,2,3\}$.
%see for detailsLemma...
%~\ref{xrxr.3},
%Lemma...
%~\ref{dr.3},
% Lemma
 %~\ref{exp(x)A3} 
% and Lemma...
 %~\ref{l.Re-Im-exp.3}.
\end{rem}

\iffalse
By \eqref{G.detC-LI.2} we have in particular
\begin{equation}
\label{Delta(f,g,h)}
\Delta(y_1,y_2)\!=\!\frac{{\rm det}
\big(I_2\!+\!\gamma(y_1,y_2)\big)}
{{\rm det}\big(I_{1}+\gamma(y_2)\big)}\!-\!1,\,\,
%\\
%&&
\Delta(y_1,y_2,y_3)\!=\!\frac{{\rm det}
\big(I_3\!+\!\gamma(y_1,y_2,y_3)\big)}
{{\rm det}\big(I_{2}+\gamma(y_2,y_3)\big)}\!-\!1.
\end{equation}
\fi

\begin{lem}[\cite{KosMor-Arx23}, Lemma 4.1]
\label{l.approx-(x,D).m=3}
If $\mu^{L_t}\perp\mu$ for all $t\in {\rm GL}(3,{\mathbb R})\setminus \{e\}
$, we can approximate at least  one of  the following  eight triplets of operators:
\begin{eqnarray*}
 &&(x_{1n},x_{2n},x_{3n}),\,\,(x_{1n},x_{2n},D_{3n}),\,\,(x_{1n},D_{2n},x_{3n}),\,\,(D_{1n},x_{2n},x_{3n}),\\
 &&(x_{1n},D_{2n},D_{3n}), (D_{1n},x_{2n},D_{3n}),\,\,(D_{1n},D_{2n},x_{3n}),\,\, (D_{1n},D_{2n},D_{3n}).
\end{eqnarray*}
\end{lem}
%\newpage
{\it We give here only some results to show how Lemma~\ref{l.det/det.m} 
%{l.m=3}
is used in the case $m=2$.}
We write compactly Lemma~\ref{xrxr.3} and  Lemma~\ref{dr.3} below as follows:
\begin{eqnarray}
&&
x_{rn}x_{rt}\,\,\eta\,\, {\mathfrak A}^3\Leftrightarrow \Delta^{(r)}=\infty,\quad
D_{rn}\,\,\eta\,\, {\mathfrak A}^3\Leftrightarrow \Delta_r=\infty,\\
\label{Delta^k(krs)}
&&\text{where}\,\,\Delta^{(r)}:=\Delta(Y_r^{(r)},Y_s^{(r)},Y_t^{(r)}),\quad
\Delta_r:=\Delta(Y_r,Y_s,Y_t),
\end{eqnarray}
and $\{r,s,t\}$ is a cyclic permutation of $\{1,2,3\}$. Here
\begin{eqnarray}
\label{Y_r^{(s)}}
&&
\Vert Y_s^{(r)} \Vert^2=\sum_{k\in \mathbb Z}
\frac{b_{rk}^2}{B_{3k}^2-(b_{1k}^2+b_{2k}^2+b_{3k}^2-b_{sk}^2)},\quad 1\leq r,s\leq 3,\\
\label{II}
&&
B_{3k}=b_{1k}+b_{2k}+b_{3k},\,\,\,\text{and}\,\,\,
\Vert Y_r\Vert^2=\sum_{k\in \mathbb Z}\frac{a_{rk}^2}{
\frac{1}{2b_{1k}}+\frac{1}{2b_{2k}}+\frac{1}{2b_{3k}}}.
\end{eqnarray}
The generators $A_{kn}$ have the following form:
\begin{equation}
\label{A(kn).m=3} 
A_{kn}=x_{1k}D_{1n}+x_{2k}D_{2n}+x_{3k}D_{3n},\,\,k,n\in \mathbb Z.
\end{equation}
\begin{lem}
\label{xrxr.3}
For any  $n,t\in {\mathbb Z}$ and $1\leq r\leq 3 $ one has
$$
x_{rn}x_{rt}{\bf 1}\in\langle A_{nk}A_{tk}{\bf 1}\mid k\in{\mathbb
Z}\rangle \,\,\Leftrightarrow\,\,
\Delta(Y_r^{(r)},Y_s^{(r)},Y_l^{(r)})=\infty,
$$
where $\{r,s,l\}$ is a cyclic permutation of $\{1,2,3\}$.
\end{lem}
For $m=3$, consider three rows as follows
\begin{equation}
 \label{lambda(D)}
\left(
\begin{array}{cccccc}
...&a_{11}&a_{12}&...&a_{1n}&...\\
...&a_{21}&a_{22}&...&a_{2n}&...\\
...&a_{31}&a_{32}&...&a_{2n}&...
\end{array}
\right)\quad
\text{and set}\quad \lambda_k=\frac{1}{2b_{1k}}+\frac{1}{2b_{2k}}+\frac{1}{2b_{3k}}.
\end{equation}
Denote by $Y_1,Y_2$ and $Y_3$ the three following vectors:
\begin{equation}
 \label{Y(1),Y(3),Y(3)}
 x_{rk}=a_{rk}/\sqrt{\lambda_k},\,\,
 %\frac{a_{rk}}{\sqrt{\lambda_k}},\,
 \,k\in{\mathbb Z}
 ,\quad Y_r=(x_{rk})_{k\in {\mathbb Z}}.
\end{equation}
\begin{lem}
\label{dr.3}
For any  $l\in {\mathbb Z}$ we have
$$
D_{rl}{\bf 1}\in\langle A_{kl}{\bf 1}\mid k\in {\mathbb Z}\rangle
\quad\Leftrightarrow\quad \Delta(Y_r,Y_s,Y_t)=\infty,
$$
where $\{r,s,t\}$ is a cyclic permutation of $\{1,2,3\}$.
\end{lem}
%\newpage
\vskip 0.3cm
\noindent{\it Acknowledgement.} 
%A.~Kosyak 
The author is very grateful to Prof. K.-H. Neeb, Prof. M.~Smirnov and 
Dr P.~Moree
for their personal efforts  to make academic  stays possible at their respective institutes. %A.~Kosyak
The author
visited: MPIM from March to April 2022 and from January to April 2023,
University of Augsburg from June to July 2022, and  University of Erlangen-Nuremberg 
from August to December 2022, all during the Russian invasion in Ukraine.
Also, Prof. R.~Kashaev  kindly invited 
%A.~Kosyak 
him to Geneva.

Further, he 
would like  to pay his respect to Prof. P. Teichner at MPIM, for his
immediate efforts started to help mathematicians
in Ukraine  after the Russian invasion.

Since the spring of 2023 A.~Kosyak is an Arnold Fellow at the
London Institute for Mathematical Sciences, and he would like to express
his gratitude to  Mrs S.~Myers Cornaby   to Miss A.~Ker Mercer and to Dr M.~Hall 
%Madeleine 
and especially to the Director of LIMS Dr T.~Fink and Prof. Y.-H.~He.
%

%%%%%%%%%%%%%%%

%%%%%%%%%%%%%%%%%
\end{document}